\documentclass[reqno]{amsart}

\usepackage{graphicx,amssymb,mathrsfs,amsmath,color,fancyhdr}
\usepackage[all]{xy}
\newtheorem{thm}{Theorem}[section]
\newtheorem{lemma}{Lemma}[section]

\newtheorem{conj}{Conjecture}[section]

\UseRawInputEncoding

\numberwithin{equation}{section}

\newcommand{\bdot}{\boldsymbol{\cdot}}

\begin{document}

\title[A reciprocity on finite abelian groups]
{A reciprocity on finite abelian groups involving zero-sum sequences II}

\begin{abstract}
Let $G$ be a finite abelian group. For any positive integers $d$ and $m$, let $\varphi_G(d)$ be the number of elements in $G$ of order $d$ and $\mathsf M(G,m)$ be the set of all zero-sum sequences of length $m$. In this paper, for any finite abelian group $H$, we prove that
$$|\mathsf M(G,|H|)|=|\mathsf M(H,|G|)|$$
if and only if $\varphi_G(d)=\varphi_H(d)$ for any $d|(|G|,|H|)$. We also consider an extension of this result to non-abelian groups in terms of invariant theory.

\end{abstract}

\author{Mao-Sheng Li}
\address{School of Mathematics, South China University of Technology, Guangzhou 510641, Guangdong, P.R. China}
\email{li.maosheng.math@gmail.com}
\author{Hanbin Zhang}
\address{School of Mathematics (Zhuhai), Sun Yat-sen University, Zhuhai 519082, Guangdong, P.R. China}
\email{zhanghb68@mail.sysu.edu.cn}

\maketitle

\section{Introduction}

Let $G$ be an additive finite abelian group of order $|G|=n$. Let $m$ be a positive integer. We call $S=g_1\bdot\ldots\bdot g_m$ a sequence over $G$ if $S$ consists of a finite sequence of terms $g_1,\ldots,g_m\in G$ which is unordered and repetition of terms is allowed. Here, $m$ is called the length of $S$ and we denote it by $|S|=m$. We define $\sigma(S)=g_1+\cdots+g_m$ and call $S$ a zero-sum sequence if $\sigma(S)$ equals 0, the identity of $G$. The studies of zero-sum sequences over finite abelian groups can be traced back to classical works of Erd\H{o}s, Ginzburg and Ziv \cite{EGZ} and Olson \cite{Olson1,Olson2}; we refer to \cite{GG} for a survey on zero-sum theory. We denote
$$\mathsf M(G,m)=\{S\text{ is a sequence over }G\text{ }|\text{ }\sigma(S)=0\text{ and }|S|=m\}.$$
Let $C_n$ be a cyclic group with $n$ elements. In 1975, Fredman \cite{Fred} proved the following very interesting reciprocity
\begin{equation}\label{Fredman}
|\mathsf M(C_n,m)|=|\mathsf M(C_m,n)|
\end{equation}
using generating functions as well as a necklace interpretation.
Later in 1999, Elashvili, Jibladze and Pataraia \cite{EJ,EJP} rediscovered the same result with different method from invariant theory. It was remarked in \cite[Introduction]{EJP} that N. Alon also independently proved (\ref{Fredman}). Meanwhile, G. Andrews, N. Alon and R. Stanley independently obtained the counting formula for $\mathsf M(C_n,m)$; see \cite[Introduction and Section 3]{EJP}. In 2011, Panyushev \cite{Pan} provided an extension of Fredman's reciprocity in terms of symmetric tensor exterior algebras.

Recall that, for any positive integer $m$, we have the following counting formula
\begin{equation}\label{eq:Formula}
    |\mathsf M(G,m)|=\frac{1}{n+m}\sum_{d|(n,m)}\varphi_G(d)\binom{n/d+m/d}{n/d},
\end{equation}
where
$$\varphi_G(d)=\sum_{\ell|d}\mu(\frac{d}{\ell})
\prod_{i=1}^r(n_i,\ell)$$ is the number of elements in $G$ of order $d$; see \cite{HZ,MRW} for proofs. It follows immediately from (\ref{eq:Formula}) that, if
\begin{equation}\label{groupcondition}
\varphi_G(d)=\varphi_H(d)\text{ for any }d| (|G|,|H|),
\end{equation}
then we have
\begin{equation}\label{generalreciprocity}
|\mathsf M(G,|H|)|=|\mathsf M(H,|G|)|.
\end{equation}
In particular, if $(|G|,|H|)=1$, then
\begin{equation}\label{eq:RationalCatalan}
|\mathsf M(G,|H|)|=|\mathsf M(H,|G|)|=\frac{1}{|G|+|H|}\binom{|G|+|H|}{|G|,|H|},
\end{equation}
which are called the rational Catalan numbers. Recently, Han and Zhang \cite{HZ} provided a bijective proof of (\ref{eq:RationalCatalan}) using rational Catalan combinatorics, which is based on a correspondence between zero-sum sequences and rational Dyck paths.

In this paper, we continue the above study and prove that (\ref{groupcondition}) is not only a sufficient condition but also a necessary condition for the reciprocity (\ref{generalreciprocity}) to hold.

\begin{thm}\label{thm:MainTheorem}
Let $G$ and $H$ be two finite abelian groups. Then we have
$$|\mathsf M(G,|H|)|=|\mathsf M(H,|G|)|$$
if and only if $\varphi_G(d)=\varphi_H(d)$ for any $d| (|G|,|H|)$.
\end{thm}

Note that, Theorem \ref{thm:MainTheorem} shows that, between any two finite abelian groups, a combinatorial property is equivalent to a group structural property.

Next, we discuss Theorem \ref{thm:MainTheorem} from the viewpoint of invariant theory and consider an extension of Theorem \ref{thm:MainTheorem} to finite groups (not necessarily abelian). Let $G$ be a finite group and $\rho: G\rightarrow GL(V)$ be a finite dimensional linear representation of $G$ over $\mathbb{C}$. Let $\mathbb{C}[V]$ denote the graded algebra of polynomial functions on $V$. We can regard $\mathbb{C}[V]$ as the symmetric algebra on $V^{*}$, the dual space of $V$. Equivalently, if $z_1,\ldots,z_n\in V^*$ is a basis, then $\mathbb{C}[V]$ is just the polynomial ring $\mathbb{C}[z_1,\ldots,z_n]$, whose elements are the homogeneous polynomials in the linear forms $z_1,\ldots,z_n$ with coefficient in $\mathbb{C}$. The action of $G$ on $V$ (through the representation $\rho$) naturally induces a right action of $G$ on $V^*$ as follows: $$x^g(v)=x(g\cdot v)=x(\rho(g)v).$$
Moreover, this action can be naturally extended to an action on $\mathbb{C}[V]$. The central topic of invariant theory is to study the algebra of polynomial invariants which is defined as follows:
$$\mathbb{C}[V]^G=\{f\in\mathbb{C}[V]\text{ }|\text{ }f^g=f,\text{ for all }g\in G\}.$$
Recall that $\mathbb{C}[V]^G$ is finitely generated and that $\mathbb{C}[V]^G$ is a graded $\mathbb{C}$-algebra; see \cite{NeuSm}. We denote by $\dim\mathbb C[V]^G_{m}$ the dimension of the $m$-th component of $\mathbb C[V]^G$ as a vector space over $\mathbb C$. In particular, if $G$ is a finite abelian group and $V$ is the regular representation of $G$ over $\mathbb{C}$, then for any positive integer $m$ we have
$$\dim \mathbb{C}[V]^G_m=|\mathsf M(G,m)|$$
(see \cite[Section 3]{HZ} for a discussion).
Now, we can restate Theorem \ref{thm:MainTheorem} in terms of invariant theory as follows.

\begin{thm}\label{Mainthm2}
Let $G$ and $H$ be two finite abelian groups. Let $V$ (resp. $V'$) be the regular representation of $G$ (resp. $H$). Then we have
$$\dim \mathbb C[V]^G_{|H|}=\dim \mathbb C[V']^H_{|G|}$$
if and only if $\varphi_G(d)=\varphi_H(d)$ for any $d| (|G|,|H|)$.
\end{thm}

For any finite group $G$ and its regular representation $V$ over $\mathbb C$, Almkvist and Fossum \cite[Section V 1.8]{AlmF} proved that
\begin{equation}\label{eq:dimFormula}
    \dim \mathbb{C}[V]^G_m=\frac{1}{n+m}\sum_{d|(n,m)}\varphi_G(d)\binom{n/d+m/d}{n/d},
\end{equation}
where
$\varphi_G(d)$ is the number of elements in $G$ of order $d$. Note that \eqref{eq:dimFormula} has the same form as \eqref{eq:Formula}. Consequently,
for any finite group $H$ and its regular representation $V'$ over $\mathbb C$, if \eqref{groupcondition} holds, then we have
\begin{equation}\label{reciprocity_finite_groups}
\dim \mathbb C[V]^G_{|H|}=\dim \mathbb C[V']^H_{|G|}.
\end{equation}
Based on this observation and Theorem \ref{thm:MainTheorem}, it is natural to consider whether \eqref{groupcondition} is a necessary condition for the reciprocity \eqref{reciprocity_finite_groups} to hold.
In fact, we prove some positive results on this problem.

\begin{thm}\label{Mainthm3}
Let $G$ and $H$ be two finite groups. Let $V$ (resp. $V'$) be the regular representation of $G$ (resp. $H$).
Assume that one of the following holds:
\begin{enumerate}
    \item $G=D_{2p}$ is the dihedral group of order $2p$, where $p\ge 5$ is a prime;
    \item $|H|\ge |G|^2$ and $|G|$ does not contain two divisors $d_1, d_2>1$ with $d_1-d_2=1$.
\end{enumerate}
Then we have
$$\dim \mathbb C[V]^G_{|H|}=\dim \mathbb C[V']^H_{|G|}$$
if and only if $\varphi_G(d)=\varphi_H(d)$ for any $d| (|G|,|H|)$.
\end{thm}

The rest of this paper is organized as follows. In Section 2, we introduce some definitions as well as some auxiliary lemmas. In Section 3, we prove our main results. In Section 4, we provide some further discussions.

\section{Preliminaries}

In this section, we introduce some definitions and notation, as well as some auxiliary results.

Let $\mathbb C$ be the field of complex numbers. Denote by $\mathbb N$ the set of positive integers and let $\mathbb N_0=\mathbb N\cup\{0\}$. Let $G$ be a finite abelian group written additively. By the fundamental theorem of finite abelian groups we have
$$G\cong C_{n_1}\oplus\cdots\oplus C_{n_r}$$
where $r=\mathsf r(G)\in \mathbb{N}_0$ is the rank of $G$, $1<n_1|\cdots|n_r\in\mathbb{N}$ are positive integers. Moreover, $n_1,\ldots,n_r$ are uniquely determined by $G$. We denote by ord$(g)$ the order of an element $g$ in a group. For any prime $p$, we denote by $Syl_p(G)$ the Sylow $p$-subgroup of $G$.

Now, we prove some auxiliary lemmas, which will be repeatedly used in the subsequent proofs. The first two technical lemmas are some inequalities involving binomial coefficients.

\begin{lemma}\label{lemma:Binom_ineq}
Let $m,n,a,b\geq 2$ be integers and $a,b|(m,n).$
\begin{enumerate}
   \item [{\rm(i)}] If $b>a$,  we have
$$
  \binom{\frac{m+n}{a}}{\frac{m}{a},\frac{n}{a}}\Big/ \binom{\frac{m+n}{b}}{\frac{m}{b},\frac{n}{b}}\geq (1+\frac{m}{n})^{n(\frac{1}{a}-\frac{1}{b})}(1+\frac{a}{b}\frac{n}{m})^{m(\frac{1}{a}-\frac{1}{b})}.
$$
Consequently, $\binom{\frac{m+n}{a}}{\frac{m}{a},\frac{n}{a}}> \binom{\frac{m+n}{b}}{\frac{m}{b},\frac{n}{b}}.$

    \item  [{\rm(ii)}] If $b\geq 2a$,  we   have $$
 a \binom{\frac{m+n}{a}}{\frac{m}{a},\frac{n}{a}}> \max\{m,n\}\binom{\frac{m+n}{b}}{\frac{m}{b},\frac{n}{b}}.
$$

\end{enumerate}

\end{lemma}
\begin{proof}
{\rm(i)}  Since $b>a$, we have

$$
 \begin{array}{rcl}
 \displaystyle\binom{\frac{m+n}{a}}{\frac{m}{a},\frac{n}{a}}&=&\displaystyle\frac{(\frac{m+n}{a})!}{\frac{m}{a}!\frac{n}{a}!}
 =\frac{
\prod_{i=\frac{m+n}{b}+1}^{ \frac{m+n}{a} }i}
{\prod_{j=\frac{m}{b}+1}^{\frac{m}{a}}j
\prod_{k=\frac{n}{b}+1}^{ \frac{n}{a}}k}
\binom{\frac{m+n}{b}}{\frac{m}{b},\frac{n}{b}}
 \\[4mm]
 &=&\displaystyle \prod_{j=\frac{m}{b}+1}^{ \frac{m}{a}  } \frac{\frac{n}{b}+j}{j}
\prod_{k=\frac{n}{b}+1}^{ \frac{n}{a}  } \frac{\frac{m}{a}+k}{k}\binom{\frac{m+n}{b}}{\frac{m}{b},\frac{n}{b}}\\[4mm]
 &\geq &\displaystyle \prod_{j=\frac{m}{b}+1}^{ \frac{m}{a}  } \frac{\frac{n}{b}+\frac{m}{a}}{\frac{m}{a}}
\prod_{k=\frac{n}{b}+1}^{ \frac{n}{a}  } \frac{\frac{m}{a}+\frac{n}{a} }{\frac{n}{a} }\binom{\frac{m+n}{b}}{\frac{m}{b},\frac{n}{b}}\\[4mm]
   &= & \displaystyle (1+\frac{m}{n})^{n(\frac{1}{a}-\frac{1}{b})}
(1+\frac{a}{b}\frac{n}{m})^{m(\frac{1}{a}-\frac{1}{b})}
\binom{\frac{m+n}{b}}{\frac{m}{b},\frac{n}{b}},
 \end{array}
 $$
as desired.

{\rm(ii)}
Without loss of generality, we assume that $m \ge n$.
Note that
 $$a\binom{\frac{m+n}{a}}{\frac{m}{a},\frac{n}{a}}=a\frac{\frac{m+n}{a} \cdots (\frac{m}{a}+2)(\frac{m}{a}+1)}{\frac{n}{a}(\frac{n}{a}-1) \cdots 2\cdot 1}\geq  (m+a)\prod_{j=2}^{ \frac{n}{a} } \frac{\frac{m}{a}+j}{j}.$$
Using the fact that $ \frac{n}{a}-1\ge \frac{n}{b}$ and that $\frac{b}{a}\geq 2\geq \frac{j+1}{j}$  (which implies $\frac{\frac{m}{a}+j+1}{j+1}\geq \frac{\frac{m}{b}+j}{j}$) for all $j\geq 1$, we have
 $$\prod_{j=2}^{ \frac{n}{a} } \frac{\frac{m}{a}+j}{j}= \prod_{j=1}^{ \frac{n}{a}-1 } \frac{\frac{m}{a}+j+1}{j+1} \geq  \prod_{j=1}^{ \frac{n}{b}  } \frac{\frac{m}{a}+j+1}{j+1}\geq  \prod_{j=1}^{ \frac{n}{b}  } \frac{\frac{m}{b}+j}{j}=\binom{\frac{m+n}{b}}{\frac{m}{b},\frac{n}{b}}.$$
 Therefore, we obtain
 $$a\binom{\frac{m+n}{a}}{\frac{m}{a},\frac{n}{a}}\geq (m+a)\prod_{j=2}^{ \frac{n}{a} } \frac{\frac{m}{a}+j}{j}> m \binom{\frac{m+n}{b}}{\frac{m}{b},\frac{n}{b}}.$$
This completes the proof.
\end{proof}

\begin{lemma}\label{lemma:Impartant_estimation}
Let $n=p^\alpha q^\beta n', m=p^\gamma q^\delta m'$ be positive integers, where $p$ and $q$ are distinct primes, $(n',pq)=(m',pq)=1$, and $\alpha, \beta, \gamma, \delta\in\mathbb N$. Suppose that $a=p^s$, $b=q^t$ ($s,t\geq 1$) satisfy $b<2a$ and  $a,b|(m,n)$. Denote
$$\Delta_{m,n}(a,b):=p^{\alpha+\gamma-2s-1}q^{\beta-t}m'n'.$$
\begin{enumerate}
   \item [{\rm(i)}]  If
  $n(\frac{1}{a}-\frac{1}{b})\geq 3$ and    $\{a,b\}\neq \{2,3\},$
  then
$$
     a\binom{\frac{m+n}{a}}{\frac{m}{a},\frac{n}{a}}
     \Big/\binom{\frac{m+n}{b}}{\frac{m}{b},\frac{n}{b}}>  2\Delta_{m,n}(a,b) q^\delta.
$$
Consequently, if $\Delta_{m,n}(a,b)\geq 1,$ then we have
\begin{equation}\label{eq:Bino_a_2b}
     a\binom{\frac{m+n}{a}}{\frac{m}{a},\frac{n}{a}}- (q^\delta-q^t) \binom{\frac{m+n}{b}}{\frac{m}{b},\frac{n}{b}}> 2b\binom{\frac{m+n}{b}}{\frac{m}{b},\frac{n}{b}}.
     \end{equation}
   Moreover, \eqref{eq:Bino_a_2b} always holds when $\{a,b\}= \{2,3\}$.
 \item [{\rm(ii)}]
  If
  $n(\frac{1}{a}-\frac{1}{b})=2,$ and   $\{a,b\}\neq \{2,3\},$
  then
    $$
 a\binom{\frac{m+n}{a}}{\frac{m}{a},\frac{n}{a}}
 \Big/\binom{\frac{m+n}{b}}{\frac{m}{b},\frac{n}{b}}>     \Delta_{m,n}(a,b) q^\delta.
$$
Consequently, if $\Delta_{m,n}(a,b)\geq 1,$ then we have
\begin{equation}\label{eq:Bino_a_b}
     a\binom{\frac{m+n}{a}}{\frac{m}{a},\frac{n}{a}}- (q^\delta-q^t) \binom{\frac{m+n}{b}}{\frac{m}{b},\frac{n}{b}}> b\binom{\frac{m+n}{b}}{\frac{m}{b},\frac{n}{b}}.
     \end{equation}
  \end{enumerate}

\end{lemma}
\begin{proof}
Without loss of generality, we assume that $m\geq n$.

(i) As $n(\frac{1}{a}-\frac{1}{b})\geq 3$, by Lemma \ref{lemma:Binom_ineq},
we have
$$
 \begin{array}{rcl}
 a\displaystyle\binom{\frac{m+n}{a}}{\frac{m}{a},\frac{n}{a}}\Big/\binom{\frac{m+n}{b}}{\frac{m}{b},\frac{n}{b}}&\geq& a (1+\frac{m}{n})^{n(\frac{1}{a}-\frac{1}{b})}(1+\frac{a}{b}\frac{n}{m})^{m(\frac{1}{a}-\frac{1}{b})}\\[4mm]
 &> & a\left(\frac{m}{n}n(\frac{1}{a}-\frac{1}{b})+\frac{m^2}{n^2}n(\frac{1}{a}-\frac{1}{b})  \right) n(\frac{1}{a}-\frac{1}{b})\frac{a}{b}\\[4mm]
 &\geq & 2amn(\frac{1}{a}-\frac{1}{b})^2\frac{a}{b}\\[4mm]
   &= & \displaystyle 2 \frac{mn}{pa^2bq^\delta } q^\delta   \left((b-a)^2 \frac{pa^2}{b^2}\right)\\[4mm]
   &= & 2 \Delta_{m,n}(a,b) q^\delta   \left((b-a)^2 \frac{pa^2}{b^2}\right).
 \end{array}
 $$
It suffices to prove that $(b-a)^2\frac{pa^2}{b^2}\geq 1$. If $b-a\geq 2$, as $b<2a$, then we have $(b-a)^2 \frac{pa^2}{b^2} \geq \frac{4a^2}{b^2}p>p$. Therefore, we may assume that $b=a+1.$ In this setting, as $\{a,b\}\neq \{2,3\}$, we have $a\geq 3$ and $(b-a)^2 \frac{pa^2}{b^2}=\frac{p a^2}{(a+1)^2} =\frac{p}{(1+1/a)^2}\geq \frac{9}{16}p>1.$

Now, we consider the case $\{a,b\}=\{2,3\}$.  In this case, we have  $(b-a)^2 \frac{pa^2}{b^2}=\frac{8}{9}$, where $p=2$. As $n(\frac{1}{2}-\frac{1}{3})\geq 3$, we have $n\geq 18$. Now, it is easy to check that $\Delta_{m,n}(a,b)\geq 2$ and therefore $\Delta_{m,n}(a,b)(b-a)^2  \frac{pa^2}{b^2}>1$. Consequently, we have
$$
 \begin{array}{rcl}
 a\displaystyle\binom{\frac{m+n}{a}}{\frac{m}{a},\frac{n}{a}}>  2 \Delta_{m,n}(a,b) q^\delta   \left((b-a)^2 \frac{pa^2}{b^2}\right) \binom{\frac{m+n}{b}}{\frac{m}{b},\frac{n}{b}} > 2  q^\delta \binom{\frac{m+n}{b}}{\frac{m}{b},\frac{n}{b}},
 \end{array}
 $$
and (\ref{eq:Bino_a_2b}) follows.

(ii) In this case, by Lemma \ref{lemma:Binom_ineq}, we have
$$
 \begin{array}{rcl}
 a\displaystyle\binom{\frac{m+n}{a}}{\frac{m}{a},\frac{n}{a}}\Big/\binom{\frac{m+n}{b}}{\frac{m}{b},\frac{n}{b}}&\geq& a (1+\frac{m}{n})^{n(\frac{1}{a}-\frac{1}{b})}(1+\frac{a}{b}\frac{n}{m})^{m(\frac{1}{a}-\frac{1}{b})}\\[4mm]
  &> & a m (\frac{1}{a}-\frac{1}{b})   n(\frac{1}{a}-\frac{1}{b})\frac{a}{b}\\[4mm]
 &= &  amn(\frac{1}{a}-\frac{1}{b})^2\frac{a}{b}\\[4mm]
 &= & \Delta_{m,n}(a,b) q^\delta   \left((b-a)^2 \frac{pa^2}{b^2}\right).
 \end{array}
 $$
Similar to the above, we have $(b-a)^2\frac{pa^2}{b^2}\geq 1$ and the desired result follows.
\end{proof}

The following three lemmas, whose proofs are based on the structure of finite abelian groups, are very useful in our subsequent proofs.

\begin{lemma}\label{lemma:min_ord}
Let $G$ be a finite abelian group of order $n$ and $H$ a finite abelian group of order $m$.  Denote
$$\begin{array}{c}
\mathcal{E}_G:=\{d\in \mathbb N\ |\ \varphi_G(d) > \varphi_H(d)\text{ for }d|(|G|,|H|)\},\\[2mm]
\mathcal{E}_H:= \{d\in \mathbb N\ |\ \varphi_G(d)< \varphi_H(d)\text{ for }d|(|G|,|H|)\}.
\end{array}$$
If $\mathcal{E}_G$ (resp. $\mathcal{E}_H$) is nonempty, then $\min{\mathcal{E}_G}$ (resp. $\min{\mathcal{E}_H}$) is a prime power.
\end{lemma}

\begin{proof} Let
$$n=\prod_{i=1}^\ell p_i^{n_i},\quad m=\prod_{i=1}^\ell p_i^{m_i},$$
where $n_i,m_i\geq 0$ for $1\le i\le \ell$.
First, we assume that $\mathcal{E}_G$ is nonempty.

For any $d|(n,m)$, let $d=\prod_{i=1}^\ell p_i^{d_i}$, where $d_i\geq 0$  and $d_i\leq \min\{m_i,n_i\}$ for $1\le i\le \ell$.
Note that,
$$G=Syl_{p_1}(G)\oplus\cdots\oplus Syl_{p_{\ell}}(G).$$
Therefore,
$g=(g_1,g_2,\cdots,g_\ell)\in G$  (where $g_i\in  Syl_{p_i}(G)$) has order $d$ if and only if ord$(g_i)=p_i^{d_i}$.
It follows that
$$\varphi_G(d)=\prod_{i=1}^\ell \varphi_{Syl_{p_i}(G)}(p_i^{n_i}).$$  Consequently,  if  $\varphi_G(d) > \varphi_H(d)$, we must have
$$\varphi_G(p_i^{d_i})= \varphi_{Syl_{p_i}(G)}(p_i^{d_i}) > \varphi_{Syl_{p_i}(H)}(p_i^{d_i})=\varphi_{H}(p_i^{d_i})$$
for some $i\in \{1,\cdots, \ell\}.$ The desired result follows immediately. It is easy to see that the proof is similar when $\mathcal{E}_H$ is nonempty.
\end{proof}

\begin{lemma}\label{lemma:min-order_p}
Let $G$ be a finite abelian group of order $n$ and $H$ a finite abelian group of order $m$. Let $n=q^\beta n'$  and $m=q^\delta m'$ with $(m',q)=(n',q)=1$. Let
$$\begin{array}{c}
\mathcal{E}=\{k\in\mathbb{N} \ \big |\ \varphi_G(q^k) \neq \varphi_H(q^k),\ q^k|(m,n) \}.
\end{array}$$
Suppose that $\mathcal{E}$ is nonempty and let $t=\min\mathcal{E}$. If $\varphi_G(q^t) < \varphi_H(q^t),$ then we have $q^{t+1}|m$, i.e., $\delta>t$.
Moreover,
$$ q^t\leq  \varphi_H(q^t)-\varphi_G(q^t)  \leq  q^\delta -q^t.$$
\end{lemma}
\begin{proof}
Let
$$Syl_q(G)= C_{q^{n_1}}\oplus\cdots\oplus C_{q^{n_c}},\quad
Syl_q(H)=C_{q^{m_1}}\oplus\cdots\oplus C_{q^{m_d}},$$
where $1\leq n_1\leq \cdots \leq n_c$ and  $1\leq m_1\leq \cdots \leq m_d.$ As  $\varphi_G(q^t) < \varphi_H(q^t)$, we have $m_d\geq t.$  We claim that $d\geq 2.$ Otherwise, we also have $c=1$.  In other words, both $Syl_q(G)$ and $Syl_q(H)$ are cyclic groups of order $\geq p^t$. So $\varphi_G(q^t)=q^t-q^{t-1}= \varphi_H(q^t)$  which contradicts our assumption. This proves the claim. Note that $d\geq 2$ implies $q^{m_1+m_d}| m,$ that is $q^{t+1}|m.$

It is easy to show that $\sum_{i=0}^t \varphi_G(q^i)=q^{e_1}$ and  $\sum_{j=0}^t \varphi_H(q^j)=q^{e_2}$ for some $e_1$ and $e_2$ with $ t\leq e_1<e_2 \leq  \delta.$  As $\varphi_H(q^i)=\varphi_G(q^i)$ for $i=0,\cdots,t-1,$ we have
$$ \varphi_H(q^t)-\varphi_G(q^t) =\sum_{i=0}^t\left(\varphi_H(q^i)-\varphi_G(q^i)\right)= \sum_{i=0}^t\varphi_H(q^i) - \sum_{j=0}^t\varphi_G(q^j)=q^{e_2}-q^{e_1}.$$
Note that $q^t\leq q^{e_1}(q^{e_2-e_1}-1)=q^{e_2}-q^{e_1}\leq q^\delta- q^t.$ Therefore, we obtain
$$ q^t\leq  \varphi_H(q^t)-\varphi_G(q^t)  \leq  q^\delta -q^t.$$
This completes the proof.
\end{proof}

\begin{lemma}\label{lemma:p_order_estimation}
Let $G$ be a finite abelian group of order $n$ and $H$ a finite abelian group of order $m$. Let $n=p^\alpha n'$  and $m=p^\gamma m'$ with $(n',p)=(m',p)=1$. Let
$$\begin{array}{c}
\mathcal{E}=\{k\in\mathbb{N} \ \big |\ \varphi_G(p^k) \neq \varphi_H(p^k),\ p^k|(m,n) \}.
\end{array}$$
Suppose that $\mathcal{E}$ is nonempty and let $s=\min\mathcal{E}$. Assume that the following hold
\begin{enumerate}
\item  $s\geq 2$;
\item  $p^{s+1}|(m,n)$;
\item  $\varphi_G(p^s) > \varphi_H(p^s)$, but  $\varphi_G(p^{s+1}) < \varphi_H(p^{s+1})$.
\end{enumerate}
Then we have $\alpha \geq s+2$ and $\gamma\geq s+2$.
\end{lemma}
\begin{proof}
Let
$$Syl_p(G)= C_{p^{n_1}}\oplus\cdots\oplus C_{p^{n_c}},\quad
Syl_p(H)=C_{p^{m_1}}\oplus\cdots\oplus C_{p^{m_d}},$$
where $1\leq n_1\leq \cdots \leq n_c$ and  $1\leq m_1\leq \cdots \leq m_d$. Similar to the proof of Lemma \ref{lemma:min-order_p}, we have $c\geq 2$ and $n_c\geq s$, as $\varphi_G(p^k)=\varphi_H(p^k)$ for $k<s$ and  $\varphi_G(p^s)>\varphi_H(p^s)$. Since $\varphi_G(p^{s+1}) < \varphi_H(p^{s+1})$,  we have  $m_d\geq s+1$. If $d=1,$ then $\varphi_H(p)=p-1.$ Therefore, we have $\varphi_G(p)= p^c-1\geq p^2-1>p-1=\varphi_H(p)$ and $s\leq 1$, which contradicts the assumption that $s\geq2$. Consequently, we have $d\geq 2$ and $p^{m_1+m_{d}}|m$, that is, $\gamma\geq s+2$.

Note that $\sum_{i=0}^s \varphi_G(p^i)=p^{e_1}$ and  $\sum_{j=0}^s \varphi_H(p^j)=p^{e_2}$ for some positive integers $e_1$ and $e_2$.  As $d\geq 2$ and $m_d\geq s+1$, we have $p^{e_2}\geq p^{s+1}$. Moreover, by the definition of $s$, we have $p^{e_1}>p^{e_2}$. Consequently, we obtain $e_1>s+1$, which implies $\alpha\geq s+2$.
\end{proof}

\section{Proofs of the main results}

In this section, we finish the proofs of our main results.

\smallskip

{\sl Proof of Theorem \ref{thm:MainTheorem}:}  It suffices to prove that if $|\mathsf M(G,|H|)|=|\mathsf M(H,|G|)|$, then we have $\varphi_G(d)=\varphi_H(d)$ holds for any $d|(|G|,|H|)$. Assume to the contrary that there exists some $d|(|G|,|H|)$ such that $\varphi_G(d)\neq\varphi_H(d)$. Let $|G|=n$ and $|H|=m$.

Let
$$a=\min\{d\ |\ \varphi_G(d)\neq \varphi_H(d)\text{ for }d|(m,n)\}.$$
Without loss of generality, we assume that
$\varphi_G(a)> \varphi_H(a).$
In this case, we shall prove that
\begin{equation}\label{largereciprocity}
|\mathsf M(G,|H|)|>|\mathsf M(H,|G|)|,
\end{equation}
which contradicts our assumption.

If $\varphi_G(d)\ge \varphi_H(d)$ holds for any $d|(m,n)$, then the desired result follows immediately. Therefore, we assume that
 $\varphi_G(d)< \varphi_H(d)$ holds for some $d|(m,n)$ and
let
$$b=\min\{d\ |\ \varphi_G(d)< \varphi_H(d)\text{ for }d|(m,n)\}.$$
 Clearly, we have $a<b$.
Moreover, by Lemma \ref{lemma:min_ord}, both $a$ and $b$ are prime powers.
Recall that
$$\mathcal{E}_H=\{e\in \mathbb{N}\  \big  |\    \varphi_G(e)< \varphi_H(e) \text{ and } e|(m,n) \}.$$
We denote
$$\mathcal{F}:=\{e\in \mathbb{N}\  \big  |\    \varphi_G(e)> \varphi_H(e),\ e>a,    \text{ and } e|(m,n) \},$$
and
$$S_{\mathcal{F}}:=\sum_{
    e\in \mathcal{F}
 } \left(\varphi_G(e)-\varphi_H(e)\right)\binom{\frac{m+n}{e}}{\frac{m}{e}, \frac{n}{e}}.$$
It is clear that $S_{\mathcal{F}}>0$. Therefore, we have
$$\begin{array}{rcl}
& &(m+n)\displaystyle\left(|\mathsf M(G,|H|)|-|\mathsf M(H,|G|)|\right) \\[4mm]
&=& \displaystyle\sum_{d| (m,n)} \left(\varphi_G(d)-\varphi_H(d)\right) \binom{\frac{m+n}{d}}{\frac{m}{d}, \frac{n}{d}}\\[4mm]
&= &  \displaystyle  \left(\varphi_G(a)-\varphi_H(a)\right)\binom{\frac{m+n}{a}}{\frac{m}{a},\frac{n}{a}} +S_{\mathcal{F}} +\sum_{
    e\in \mathcal{E}_H
 } \left(\varphi_G(e)-\varphi_H(e)\right)\binom{\frac{m+n}{e}}{\frac{m}{e}, \frac{n}{e}} \\[4mm]
 &= & \displaystyle  \left(\varphi_G(a)-\varphi_H(a)\right)
 \binom{\frac{m+n}{a}}{\frac{m}{a},\frac{n}{a}}+S_{\mathcal{F}} -\sum_{
    e\in \mathcal{E}_H
 }(\varphi_H(e)- \varphi_G(e))\binom{\frac{m+n}{e}}{\frac{m}{e}, \frac{n}{e}} \\[4mm]
  &\geq & \displaystyle a\binom{\frac{m+n}{a}}{\frac{m}{a},\frac{n}{a}}+S_{\mathcal{F}}-\sum_{
    e\in \mathcal{E}_H
 } (\varphi_H(e)- \varphi_G(e))\binom{\frac{m+n}{e}}{\frac{m}{e}, \frac{n}{e}},\\
\end{array}
$$
where the last inequality follows from Lemma \ref{lemma:min-order_p}.

Therefore, in order to prove (\ref{largereciprocity}), it suffices to show that
\begin{equation}\label{sufficestoprove}
a\binom{\frac{m+n}{a}}{\frac{m}{a},\frac{n}{a}}+S_{\mathcal{F}}  > \sum_{
    e\in \mathcal{E}_H
 } (\varphi_H(e)- \varphi_G(e))\binom{\frac{m+n}{e}}{\frac{m}{e}, \frac{n}{e}}.
\end{equation}
As $a$ and $b$ are prime powers, we may assume that $a=p^s$ and $b=q^t$ ($p$ and $q$ are primes, but not necessarily distinct).

If $b\ge 2a$, i.e., $q^t\geq 2 p^s$. By Lemma \ref{lemma:Binom_ineq}.(ii), we have
$$a\binom{\frac{m+n}{a}}{\frac{m}{a},\frac{n}{a}}  > m \binom{\frac{m+n}{b}}{\frac{m}{b}, \frac{n}{b}}\geq \sum_{
    e\in \mathcal{E}_H
 } \varphi_H(e)\binom{\frac{m+n}{b}}{\frac{m}{b}, \frac{n}{b}}\geq \sum_{
    e\in \mathcal{E}_H
 } (\varphi_H(e)- \varphi_G(e))\binom{\frac{m+n}{e}}{\frac{m}{e}, \frac{n}{e}},$$
and (\ref{sufficestoprove}) follows immediately.

If $a<b< 2a$, i.e., $p^s<q^t< 2 p^s$. In this case, we have $$q^{t-1}=\frac{q^t}{q}< \frac{2p^s}{q}\leq p^s.$$
Consequently, we have $p\neq q$. Moreover, by the definition of $a=p^s$, we have $t=\min\{ i\in\mathbb N \ \big | \ \varphi_G(q^i)\neq  \varphi_H(q^i)\}$ and  $\varphi_G(q^t)<  \varphi_H(q^t)$. By Lemma \ref{lemma:min-order_p}, we have $q^{t+1}|m$ and $p^{s+1}|n$.  Assume that
$$n=p^\alpha q^\beta n',\quad m=p^\gamma q^\delta m',$$
where $(n',pq)=(m',pq)=1$.
Then we have  $a=p^s<b=q^t<2p^s$  and  $1\leq s< \alpha,$ $1\leq t< \delta$, and  $s\leq \gamma$, $t\leq \beta$. Consequently, we have $p|n(\frac{1}{a}-\frac{1}{b})$, which implies $n(\frac{1}{a}-\frac{1}{b})\geq p\geq 2$. We distinguish two cases.

{\bf Case 1:} Assume that  $n(\frac{1}{a}-\frac{1}{b})=2$, that is, $b-a=1$ and $p=2$. In this case, we have $n=p^{s+1}q^t.$

{\bf Subcase 1.1:} Assume that $s\geq 2$. First, we claim that $\gamma=s.$ In fact, as $p^s|(m,n)$, we have $\gamma\geq s$.  If $\gamma\geq s+1,$ then $$\sum_{i=0}^{s+1}\varphi_H({p^i})\geq p^{s+1}=\sum_{i=0}^{s+1}\varphi_G({p^i})>\sum_{i=0}^{s}\varphi_H({p^i}).$$   Therefore, $\varphi_H({p^{s+1}})>0.$ On the other hand, we have $\varphi_G({p^{s+1}})=0$, as $Syl_p(G)$ can not be a cyclic group. Hence, $\varphi_H({p^{s+1}})>\varphi_G({p^{s+1}}).$ By Lemma \ref{lemma:p_order_estimation}, we obtain $p^{s+2}|n$, a contradiction. So $\gamma\leq s.$ This completes the proof of the claim.

Let $d|(m,n).$ In this case,  $\varphi_G(d)>\varphi_H(d)$  if and only if $d=p^sq^i$ for $i=0,1,\cdots, t,$ and  $\varphi_G(d)<\varphi_H(d)$  if and only if $d=p^jq^t$ for $j=0,1,\cdots, s.$ By Lemma \ref{lemma:min-order_p}, we have  $\varphi_H(b)- \varphi_G(b)\leq q^\delta -q^{t}.$
Recall that $\Delta_{m,n}(a,b)=p^{\alpha+\gamma-2s-1}q^{\beta-t}m'n'$. It is easy to see that $\Delta_{m,n}(a,b)\geq 1$. Note that $\{a,b\}\neq \{2,3\}$, as $a=p^s\geq 4$. Since $n(\frac{1}{a}-\frac{1}{b})=p= 2$, by Lemma \ref{lemma:Impartant_estimation}.(ii), we have
$$\begin{array}{rcl}
& &\displaystyle a\binom{\frac{m+n}{a}}{\frac{m}{a},\frac{n}{a}} - \sum_{
    e\in \mathcal{E}_H
 } (\varphi_H(e)- \varphi_G(e))\binom{\frac{m+n}{e}}{\frac{m}{e}, \frac{n}{e}}\\[4mm]
    &\geq & \displaystyle   a\binom{\frac{m+n}{a}}{\frac{m}{a},\frac{n}{a}}-  (q^\delta-q^{t}) \binom{\frac{m+n}{b}}{\frac{m}{b},\frac{n}{b}}- \sum_{
    j=1}^t \varphi_H(p^jb)\binom{\frac{m+n}{p^jb}}{\frac{m}{p^jb},\frac{n}{p^jb}} \\[4mm]
    &\geq &\displaystyle    b \binom{\frac{m+n}{b}}{\frac{m}{b},\frac{n}{b}}- \sum_{
    j=1}^t \varphi_H(p^jb)\binom{\frac{m+n}{pb}}{\frac{m}{pb},\frac{n}{pb}}\\[4mm]
    &\geq &\displaystyle   b \binom{\frac{m+n}{b}}{\frac{m}{b},\frac{n}{b}}- m\binom{\frac{m+n}{pb}}{\frac{m}{pb},\frac{n}{pb}}>0,
\end{array}
$$
where the last inequality follows from Lemma \ref{lemma:Binom_ineq}.(ii). Therefore, (\ref{sufficestoprove}) follows.

{\bf Subcase 1.2:}
Assume that $s=1$. Therefore, $p^s=2$ and $2=p^s<q^t<2p^s=4$ which implies $q^t=3$  and $n=p^{s+1}q^t=12$. In this case, $G=C_2\oplus C_6$ and $H=C_{2^\gamma}\oplus H'$, where $|H'|=3^\delta m'$.
It is easy to verify that
$
 2\binom{\frac{m}{2}+6}{6} > m \binom{\frac{m}{3}+4}{4}$ for all $m\geq 18$ and $6|m$. Therefore, we have
$$\begin{array}{rcl}
& & \displaystyle  a\binom{\frac{m+n}{a}}{\frac{m}{a},\frac{n}{a}} - \sum_{
    e\in \mathcal{E}_H
 } (\varphi_H(e)- \varphi_G(e))\binom{\frac{m+n}{e}}{\frac{m}{e}, \frac{n}{e}}\\[4mm]
    &\geq  & \displaystyle  2 \binom{  \frac{m+n}{2}}{\frac{n}{2}}  - \sum_{e\in \mathcal{E}_H}   \varphi_H(e) \binom{  \frac{m+n}{3}}{\frac{n}{3}} \\[4mm]
    &\geq  & \displaystyle  2 \binom{  \frac{m}{2}+6}{6}  - m   \binom{\frac{m}{3}+4}{4}>0.
\end{array}
$$
Consequently, (\ref{sufficestoprove}) follows.

{\bf Case 2:}
Assume that $n(\frac{1}{a}-\frac{1}{b})\geq 3$. Recall that $\Delta_{m,n}(a,b)=p^{\alpha+\gamma-2s-1}q^{\beta-t}m'n'$.  As $\alpha+\gamma>2s,$ $\Delta_{m,n}(a,b)\geq 1.$ By Lemma \ref{lemma:Impartant_estimation}.(i),
 we have
\begin{equation}\label{Case2_a2b}
 a\binom{\frac{m+n}{a}}{\frac{m}{a},\frac{n}{a}}- (q^\delta-q^{t} )  \binom{\frac{m+n}{b}}{\frac{m}{b},\frac{n}{b}}>  2b  \binom{\frac{m+n}{b}}{\frac{m}{b},\frac{n}{b}}.
\end{equation}
We denote
$$\mathcal{E}_1:=\{e\in\mathcal{E}_H\ \big | \ b<e<2b\},\quad  \mathcal{E}_2:=\{e\in\mathcal{E}_H\ \big | \ e\geq 2b\}.$$ By Lemma \ref{lemma:min-order_p}, we have  $\varphi_H(b)- \varphi_G(b)\leq q^\delta -q^{t}.$ Therefore, by (\ref{Case2_a2b}), we obtain
\begin{align*}
   &\quad\ S_{\mathcal F}+\displaystyle a\binom{\frac{m+n}{a}}{\frac{m}{a},\frac{m}{a}} - \sum_{
    e\in \mathcal{E}_H
 } (\varphi_H(e)- \varphi_G(e))\binom{\frac{m+n}{e}}{\frac{m}{e}, \frac{n}{e}} \\
 &\geq S_{\mathcal F}+\displaystyle a\binom{\frac{m+n}{a}}{\frac{m}{a},\frac{m}{a}} - (\varphi_H(b)- \varphi_G(b))\binom{\frac{m+n}{b}}{\frac{m}{b}, \frac{n}{b}} -\sum_{
    e\in \mathcal{E}_1\cup  \mathcal{E}_2
 } \varphi_H(e)\binom{\frac{m+n}{e}}{\frac{m}{e}, \frac{n}{e}}\\
  &>  \displaystyle S_{\mathcal F}+ 2b \binom{\frac{m+n}{b}}{\frac{m}{b}, \frac{n}{b}} -\sum_{
    e_1\in \mathcal{E}_1
 } \varphi_H(e_1)\binom{\frac{m+n}{e_1}}{\frac{m}{e_1}, \frac{n}{e_1}}-\sum_{
    e_2\in \mathcal{E}_2
 } \varphi_H(e_2)\binom{\frac{m+n}{e_2}}{\frac{m}{e_2}, \frac{n}{e_2}}\\
   &=  \displaystyle S_{\mathcal F}+b \binom{\frac{m+n}{b}}{\frac{m}{b}, \frac{n}{b}} -\sum_{
    e_1\in \mathcal{E}_1
 } \varphi_H(e_1)\binom{\frac{m+n}{e_1}}{\frac{m}{e_1}, \frac{n}{e_1}}+b \binom{\frac{m+n}{b}}{\frac{m}{b}, \frac{n}{b}} -\sum_{
    e_2\in \mathcal{E}_2
 } \varphi_H(e_2)\binom{\frac{m+n}{e_2}}{\frac{m}{e_2}, \frac{n}{e_2}}.
\end{align*}
Denote
$$\mathcal S_1:=S_{\mathcal F}+b \binom{\frac{m+n}{b}}{\frac{m}{b}, \frac{n}{b}} -\sum_{
    e_1\in \mathcal{E}_1
 } \varphi_H(e_1)\binom{\frac{m+n}{e_1}}{\frac{m}{e_1}, \frac{n}{e_1}},$$
and
$$
 \mathcal S_2:=b \binom{\frac{m+n}{b}}{\frac{m}{b} \frac{n}{b}} -\sum_{
    e_2\in \mathcal{E}_2
 } \varphi_H(e_2)\binom{\frac{m+n}{e_2}}{\frac{m}{e_2}, \frac{n}{e_2}}.$$
In order to prove (\ref{largereciprocity}), it suffices to show that $\mathcal S_1>0$ and $\mathcal S_2>0$.

First, we consider $\mathcal S_2$, as it is easier to handle. If $\mathcal{E}_2$  is empty, then the desired result follows. Assume that $\mathcal{E}_2$  is not empty, let $c=\min \mathcal{E}_2.$ By definition, we have $c\geq 2b$. Therefore, by Lemma \ref{lemma:Binom_ineq}.(ii), we have
$$
\begin{array}{rcl}
\displaystyle b \binom{\frac{m+n}{b}}{\frac{m}{b}, \frac{n}{b}} -\sum_{
    e_2\in \mathcal{E}_2
 } \varphi_H(e_2)\binom{\frac{m+n}{e_2}}{\frac{m}{e_2}, \frac{n}{e_2}} &\geq& \displaystyle  b \binom{\frac{m+n}{b}}{\frac{m}{b}, \frac{n}{b}} -\sum_{
    e_2\in \mathcal{E}_2
 } \varphi_H(e_2)\binom{\frac{m+n}{c}}{\frac{m}{c}, \frac{n}{c}}\\[5mm]
 &\geq & \displaystyle     b \binom{\frac{m+n}{b}}{\frac{m}{b}, \frac{n}{b}} - m \binom{\frac{m+n}{c}}{\frac{m}{c}, \frac{n}{c}}>0.
\end{array}
$$
Therefore, $\mathcal S_2>0$, as desired.

Next, we consider $\mathcal S_1$. If $\mathcal{E}_1$ is empty, then the desired result follows. So, we may assume that $\mathcal{E}_1$ is not empty. First, we claim that
$\mathcal{E}_1$ consists of powers of distinct primes, i.e., $\mathcal{E}_1=\{q_1^{t_1},q_2^{t_2},\cdots,q_v^{t_v}\}$, where $q_i\neq q_j$ for $i\neq j$. For each $e_1\in\mathcal{E}_1,$ by definition,  $\varphi_G(e_1)< \varphi_H(e_1)$ and $q^t<e_1<2q^t$. Therefore, there is some prime power $\ell^k$ such that $e_1=\ell^k e_1'$ ($e_1'\in \mathbb{N}$) and $\varphi_G(\ell^k)< \varphi_H(\ell^k).$  If $e_1'\geq 2,$ we have $\ell^k=\frac{e_1}{2}<\frac{2q^t}{2}=q^t$ which contradicts the definition of $b=q^t$. Therefore, we have $e_1'=1$. This completes the proof of the claim.

Denote $$\mathcal{E}_1^{G}:=\{\ell^k\in \mathcal{E}_1  \ \big |\ \ell\neq p\   \text{and }\varphi_G(\ell^i)> \varphi_H(\ell^i)\text{ for some }i<k \}.$$

For any $\ell^k\in\mathcal{E}_1^{G}$, let $u=\min\{i\in\mathbb N\ |\ \varphi_G(\ell^i)> \varphi_H(\ell^i)\}$. Moreover, by the definitions of $b$ and $\mathcal{E}_1$, we also have
$$u=\min\{i\in\mathbb N\ |\ \varphi_G(\ell^i)\neq \varphi_H(\ell^i)\}.$$
By Lemma \ref{lemma:min-order_p}, we have $\varphi_G(\ell^u)-\varphi_H(\ell^u)\ge \ell^u$. By Lemma \ref{lemma:Binom_ineq}.(ii), we have
\begin{equation}\label{sf_larger_than_e1}
(\varphi_G(\ell^u)-\varphi_H(\ell^u))\binom{\frac{m+n}{\ell^u}}{\frac{m}{\ell^u}, \frac{n}{\ell^u}}
\ge \ell^u\binom{\frac{m+n}{\ell^u}}{\frac{m}{\ell^u}, \frac{n}{\ell^u}}
>\varphi_H(\ell^k)\binom{\frac{m+n}{\ell^k}}{\frac{m}{\ell^k}, \frac{n}{\ell^k}}.
\end{equation}
Recall that
$S_{\mathcal{F}}=\sum_{
    e\in \mathcal{F}
 } \left(\varphi_G(e)-\varphi_H(e)\right)\binom{\frac{m+n}{e}}{\frac{m}{e}, \frac{n}{e}}$.
Since $\mathcal{E}_1^{G}$ consists of powers of distinct primes (which are different from $p,q$), by (\ref{sf_larger_than_e1}), we obtain
$$
\displaystyle \begin{array}{rcl}
&&\displaystyle S_{\mathcal F}+b \binom{\frac{m+n}{b}}{\frac{m}{b}, \frac{n}{b}} -\sum_{
    e_1\in \mathcal{E}_1
 } \varphi_H(e_1)\binom{\frac{m+n}{e_1}}{\frac{m}{e_1}, \frac{n}{e_1}}\\[4mm]
 &= &\displaystyle S_{\mathcal F}-\sum_{
    e\in \mathcal{E}_1^{G}
 } \varphi_H(e)\binom{\frac{m+n}{e}}{\frac{m}{e}, \frac{n}{e}}+\displaystyle b \binom{\frac{m+n}{b}}{\frac{m}{b}, \frac{n}{b}} -\sum_{
    e_1\in \mathcal{E}_1\setminus \mathcal{E}_1^{G}
 } \varphi_H(e_1)\binom{\frac{m+n}{e_1}}{\frac{m}{e_1}, \frac{n}{e_1}}\\[5mm]
 &\geq &\displaystyle b \binom{\frac{m+n}{b}}{\frac{m}{b}, \frac{n}{b}} -\sum_{
    e_1\in \mathcal{E}_1\setminus \mathcal{E}_1^{G}
 } \varphi_H(e_1)\binom{\frac{m+n}{e_1}}{\frac{m}{e_1}, \frac{n}{e_1}} .
 \end{array}
 $$
As a result, in order to prove $\mathcal S_1>0$, it suffices to show that
\begin{equation}\label{sufficestoprove_S2}
b \binom{\frac{m+n}{b}}{\frac{m}{b}, \frac{n}{b}} -\sum_{
    e_1\in \mathcal{E}_1\setminus \mathcal{E}_1^{G}
 } \varphi_H(e_1)\binom{\frac{m+n}{e_1}}{\frac{m}{e_1}, \frac{n}{e_1}}>0.
\end{equation}

If $\mathcal{E}_{1}\setminus \mathcal{E}_1^{G}$ is empty, then clearly we have $\mathcal S_1>0$. Therefore, suppose that $\mathcal{E}_{1}\setminus \mathcal{E}_1^{G}$ is nonempty. Without loss of generality, we may assume that
$$\mathcal{E}_{1}\setminus \mathcal{E}_1^{G}=\{q_1^{t_1},\cdots,q_L^{t_L}\},$$
where $L\le v$ and $q_1^{t_1}<\cdots< q_L^{t_L}$. We denote $b_0:=b$ (with $q_0:=q$ and $t_0:=t$) and $b_i:=q_i^{t_i}$ for $i=1,2,\cdots,L$. Therefore, $b_0<b_1<b_2<\cdots<b_L<2b=2q^t$. For each $i=0,1,\cdots, L-1,$ let
$$n=q_i^{\beta_i}q_{i+1}^{\beta_{i+1}} n_i,\quad m=q_i^{\delta_i}q_{i+1}^{\delta_{i+1}} m_i,$$ where $(n_i,q_iq_{i+1})=(m_i,q_iq_{i+1})=1,$ and denote
$$\Delta_{m,n}(b_i,b_{i+1}):= q_i^{\beta_i+\delta_{i}-2t_i-1}q_{i+1}^{\beta_{i+1}-t_{i+1}}m_in_i.$$

In the following, we shall prove that
\begin{equation}\label{eq:summand_term_i}
    \displaystyle b_i \binom{\frac{m+n}{b_i}}{\frac{m}{b_i}, \frac{n}{b_i}} -\varphi_H(b_{i+1})\binom{\frac{m+n}{b_{i+1}}}{\frac{m}{b_{i+1}}, \frac{n}{b_{i+1}}}> b_{i+1}\binom{\frac{m+n}{b_{i+1}}}{\frac{m}{b_{i+1}}, \frac{n}{b_{i+1}}},\quad i=0, 1,\cdots,L-1.
\end{equation}
Note that, if (\ref{eq:summand_term_i}) holds, then we have
$$\sum_{i=0}^{L-1}\left(\displaystyle b_i \binom{\frac{m+n}{b_i}}{\frac{m}{b_i}, \frac{n}{b_i}} -\varphi_H(b_{i+1})\binom{\frac{m+n}{b_{i+1}}}{\frac{m}{b_{i+1}}, \frac{n}{b_{i+1}}}\right)>\sum_{i=0}^{L-1} b_{i+1}\binom{\frac{m+n}{b_{i+1}}}{\frac{m}{b_{i+1}}, \frac{n}{b_{i+1}}},
$$
which implies
$$b_0 \binom{\frac{m+n}{b_0}}{\frac{m}{b_0}, \frac{n}{b_0}}-\sum_{i=0}^{L-1} \varphi_H(b_{i+1})\binom{\frac{m+n}{b_{i+1}}}{\frac{m}{b_{i+1}}, \frac{n}{b_{i+1}}}>b_L \binom{\frac{m+n}{b_L}}{\frac{m}{b_L}, \frac{n}{b_L}}.$$
Therefore, we have
$$b \binom{\frac{m+n}{b}}{\frac{m}{b}, \frac{n}{b}} -\sum_{
    e_1\in \mathcal{E}_1\setminus \mathcal{E}_1^{G}
 } \varphi_H(e_1)\binom{\frac{m+n}{e_1}}{\frac{m}{e_1}, \frac{n}{e_1}} >b_L \binom{\frac{m+n}{b_L}}{\frac{m}{b_L}, \frac{n}{b_L}}>0,$$
and (\ref{sufficestoprove_S2}) follows.

Now, we prove (\ref{eq:summand_term_i}). Note that, for each $i=0,1,\cdots,L-1$,
\begin{itemize}

\item if $q_i=p$, then we have  $q_i^{t_i}=p^{s+1}$, as $q_i^{t_i}<2b<4p^s$;

\item if $q_i\neq p$, by Lemma \ref{lemma:min-order_p} (note that $t_i=\min\{j\ \big |\ \varphi_G(q_i^j)\neq \varphi_H(q_i^j)\}$ and that $\varphi_G(q_i^{t_i})< \varphi_H(q_i^{t_i})$ as $b_i=q_i^{t_i}\in \mathcal{E}_1\setminus \mathcal{E}_1^G$),  we have $\delta_i>t_i$;

\item we have $\{b_i,b_{i+1}\}\neq \{2,3\}$, as $b_i>q^t\geq 2$.

\end{itemize}
We distinguish three cases.

{\bf Subcase 2.1:} Assume that $b_{i+1}\neq p^{s+1}$.  It is easy to verify that, if  $n$ has at least three different prime divisors, then $n(\frac{1}{b_i}-\frac{1}{b_{i+1}})\geq 3$.

If $b_i\neq p^{s+1}$, as mentioned above, we have $\delta_i>t_i$. Therefore, $p$, $q_i$, and $q_{i+1}$ are different prime divisors of $n$.

If $b_i= p^{s+1},$  we have $q$, $p$, and $q_{i+1}$ are   different prime divisors  of $n.$ Moreover, we have $q^t|n_i$ and  $\frac{n_i}{q_i}\geq \frac{q^t}{q_i} \geq \frac{q^t}{p^s}>1$, where $q_i=p$.

Therefore, we always have that $n(\frac{1}{b_i}-\frac{1}{b_{i+1}})\geq 3$ and  $\Delta_{m,n}(b_i,b_{i+1})\geq 1$. By Lemma \ref{lemma:Impartant_estimation}.(i), we have
    $$\begin{array}{rcl}
     \displaystyle b_i \binom{\frac{m+n}{b_i}}{\frac{m}{b_i}, \frac{n}{b_i}}    >2 \Delta_{m,n}(b_{i},b_{i+1}) q_{i+1}^{\delta_{i+1}}\binom{\frac{m+n}{b_{i+1}}}{\frac{m}{b_{i+1}}, \frac{n}{b_{i+1}}}\geq 2 q_{i+1}^{\delta_{i+1}}\binom{\frac{m+n}{b_{i+1}}}{\frac{m}{b_{i+1}}, \frac{n}{b_{i+1}}}.
\end{array}
$$
As $b_{i+1}< q_{i+1}^{\delta_{i+1}}$ and $\varphi_H(b_{i+1})< q_{i+1}^{\delta_{i+1}}$, \eqref{eq:summand_term_i} follows.

{\bf Subcase 2.2:}  Assume that $b_{i+1}=p^{s+1}$ and $s\geq 2$. In this case, by Lemma \ref{lemma:p_order_estimation}, we have $p^{s+2}|n$. Consequently, $p|n(\frac{1}{b_{i}}-\frac{1}{b_{i+1}})$, which implies  $n(\frac{1}{b_{i}}-\frac{1}{b_{i+1}})\geq 2$.  Moreover, as $\delta_i>t_i$ and  $\beta_{i+1}\geq s+2>s+1=t_{i+1}$, we have $\Delta_{m,n}(b_i,b_{i+1})\geq p\geq 2.$
By Lemma \ref{lemma:Impartant_estimation}, we obtain
    $$\begin{array}{rcl}
     \displaystyle b_i \binom{\frac{m+n}{b_i}}{\frac{m}{b_i}, \frac{n}{b_i}}    >\Delta_{m,n}(b_{i},b_{i+1}) p^{\delta_{i+1}}\binom{\frac{m+n}{b_{i+1}}}{\frac{m}{b_{i+1}}, \frac{n}{b_{i+1}}}\geq 2 p^{\delta_{i+1}}\binom{\frac{m+n}{b_{i+1}}}{\frac{m}{b_{i+1}}, \frac{n}{b_{i+1}}}.
\end{array}
$$
As $b_{i+1}< p^{\delta_{i+1}}$ and $\varphi_H(b_{i+1})< p^{\delta_{i+1}}$, \eqref{eq:summand_term_i} follows.

{\bf Subcase 2.3:} Assume that $b_{i+1}=p^{s+1}$ and $s=1$. In this case, we have $p\leq 3$. In fact, as $b<2p^s \le p^{s+1}<2b$, we have $p/2=p^{s+1}/2p^s<2b/b=2$.

{\bf Subcase 2.3.1:} Assume that $p=2$. In this case, we have that $a=p^s=2$ and $b=q^t=3$ and $\mathcal{E}_{1}\setminus \mathcal{E}_1^{G}\subseteq \{4,5\}$. Therefore, it suffices to consider the case  $b_{0}=q^{t_0}=3$ and $b_{1}=p^{s+1}=4.$ Note that $n(\frac{1}{a}-\frac{1}{b})\geq 3.$ Therefore, $n\geq 18$.  Hence $n(\frac{1}{b_0}-\frac{1}{b_{1}})\geq \frac{18}{16}>1$ which implies that $n(\frac{1}{b_0}-\frac{1}{b_{1}})\geq 2$. If  $n(\frac{1}{3}-\frac{1}{4})=2$ (i.e., $n=24$), then $\Delta_{m,n}(b_{0},b_{1})\geq 2$ (as $ \delta_0>t_0$ and $\beta_{1}=3$ and $t_{1}=2$). Therefore, by Lemma \ref{lemma:Impartant_estimation}.(ii), we have
    $$\begin{array}{rcl}
     \displaystyle b_0 \binom{\frac{m+n}{b_0}}{\frac{m}{b_0}, \frac{n}{b_0}}    >\Delta_{m,n}(b_{0},b_{1}) 2^{\delta_{1}}\binom{\frac{m+n}{b_{1}}}{\frac{m}{b_{1}}, \frac{n}{b_{1}}}\geq 2\cdot 2^{\delta_{1}}\binom{\frac{m+n}{b_{1}}}{\frac{m}{b_{1}}, \frac{n}{b_{1}}}.
\end{array}
$$
If $n(\frac{1}{3}-\frac{1}{4})\geq 3$, as $\Delta_{m,n}(b_0,b_{1})\geq 1$, by Lemma \ref{lemma:Impartant_estimation}.(i),
    $$\begin{array}{rcl}
     \displaystyle b_{0} \binom{\frac{m+n}{b_0}}{\frac{m}{b_0}, \frac{n}{b_0}}    >2\Delta_{m,n}(b_{0},b_{1}) 2^{\delta_{1}}\binom{\frac{m+n}{b_{1}}}{\frac{m}{b_{1}}, \frac{n}{b_{1}}}\geq 2\cdot 2^{\delta_{1}}\binom{\frac{m+n}{b_{1}}}{\frac{m}{b_{1}}, \frac{n}{b_{1}}}.
\end{array}
$$
As  $b_{1}< 2^{\delta_{1}}$ and $\varphi_H(b_{1})< 2^{\delta_{1}}$, \eqref{eq:summand_term_i} follows.


{\bf Subcase 2.3.2:} Assume that $p=3$. In this case, $q^t\in\{4,5\}$.  If $q^t=4$, then $\mathcal{E}_{1}\setminus \mathcal{E}_1^G\subseteq \{5,7\}$, which contains no powers of $p$. Therefore, we only need to consider $q^t=5.$ As a result, we have $a=p^s=3$ and $b=q^t=5$, and $\mathcal{E}_{1}\setminus \mathcal{E}_1^G\subseteq \{ 7,8,9\}.$ It suffices to consider the following three cases
$$
(1):\ \{b_0,b_{1}\}=\{5,9\},\quad (2):\ \{b_0,b_{1}\}=\{7,9\},\quad (3):\ \{b_0,b_{1}\}=\{8,9\}.$$
Clearly, we have $\Delta_{m,n}(b_0,b_{1})\geq 1$, as $\delta_0>t_0$. For the case (1), we have  $n(\frac{1}{5}-\frac{1}{9})\geq 4$ as $5|n$ and $9|n$. For the cases (2) and (3), $n$ has at least three different prime divisors. Therefore, we always have $n(\frac{1}{b_0}-\frac{1}{b_{1}})\geq 3$. By Lemma \ref{lemma:Impartant_estimation}.(i),
$$\begin{array}{rcl}
     \displaystyle b_{0} \binom{\frac{m+n}{b_0}}{\frac{m}{b_0}, \frac{n}{b_0}}    >2\Delta_{m,n}(b_{0},b_{1}) 3^{\delta_{1}}\binom{\frac{m+n}{b_{1}}}{\frac{m}{b_{1}}, \frac{n}{b_{1}}}\geq 2\cdot 3^{\delta_{1}}\binom{\frac{m+n}{b_{1}}}{\frac{m}{b_{1}}, \frac{n}{b_{1}}}.
\end{array}
$$
As  $b_{1}< 3^{\delta_{1}}$ and $\varphi_H(b_{1})< 3^{\delta_{1}}$, \eqref{eq:summand_term_i} follows. This completes the proof.
\qed

\bigskip

Next, we prove Theorem \ref{Mainthm3}. Let $p\ge 5$ be a prime. Recall that the dihedral group of order $2p$ is defined as
$$D_{2p}=\langle x,y\ |\ x^2=1=y^p,\ xyx=y^{-1}\rangle.$$

{\sl Proof of Theorem \ref{Mainthm3}:}  (1) Firstly, we assume that $G=D_{2p}$ is the dihedral group of order $2p$, where $p\ge 5$ is a prime.

It suffices to prove that if $\dim \mathbb C[V]^G_{|H|}=\dim \mathbb C[V']^H_{|G|}$, then we have $\varphi_G(d)=\varphi_H(d)$ holds for any $d|(|G|,|H|)$.  As $|G|=2p$, $d|(|G|,|H|)$ if and only if $d\in\{1,2,p,2p\}$. The desired result follows easily if $(|G|,|H|)\in\{1,2,p\}$. Therefore, it suffices to consider the case  $(|G|,|H|)=2p$, i.e., $2p|m$.

Note that $\varphi_H(p)\geq \varphi_G(p)=p-1$, and $\varphi_H(2p)\geq \varphi_G(2p)=0$.
If $\varphi_H(2)\ge \varphi_G(2),$ by the formula \eqref{eq:dimFormula}, we have the desired result.
Now, suppose that $\varphi_H(2)<\varphi_G(2)=p$.  Let  $|H|=m=2\ell p$ where $\ell\in \mathbb N$. In this case, we have
$$
\begin{array}{rcl}
&& (m+2p)(\dim \mathbb C[V]^G_{|H|}-\dim \mathbb C[V']^H_{|G|})\\[3mm]
&=&   \displaystyle\sum_{d| 2p} \left(\varphi_G(d)-\varphi_H(d)\right) \binom{\frac{m+2p}{d}}{\frac{m}{d}, \frac{2p}{d}}\\[6mm]
&\geq &   \displaystyle  \left(\varphi_G(2)-\varphi_H(2)\right) \binom{\frac{m+2p}{2}}{\frac{m}{2},\frac{2p}{2}} - \varphi_H(p)  \binom{\frac{m+2p}{p}}{\frac{m}{p}, \frac{2p}{p}}-\varphi_H(2p)  \binom{\frac{m+2p}{2p}}{\frac{m}{2p}, \frac{2p}{2p}} \\[6mm]
&\geq  &   \displaystyle    \binom{\frac{m+2p}{2}}{\frac{m}{2},\frac{2p}{2}} - m  \binom{\frac{m+2p}{p}}{\frac{m}{p}, \frac{2p}{p}}=  \binom{\ell p+p}{p} - m  \binom{2\ell+2}{2}.
\end{array}
$$
Moreover, as $p\geq 5$, we have
 $$
 \begin{array}{rcl}
\displaystyle \binom{\ell p+p}{p}&=& \displaystyle\prod_{i=1}^p \frac{\ell p+i}{i}\geq \ell p (\prod_{i=2}^p \frac{\ell p+i}{i}) \\[3mm]
&>&   \displaystyle\ell p \prod_{i=2}^p (\ell+1) \geq  m \frac{(\ell+1)^4}{2}>  m  \binom{2\ell+2}{2}.
\end{array}
$$
Therefore, we obtain $\dim \mathbb C[V]^G_{|H|}>\dim \mathbb C[V']^H_{|G|},$ a contradiction.

\bigskip

(2) Secondly, we assume that $|H|\ge |G|^2$ and $|G|$ does not contain two divisors $d_1, d_2>1$ with $d_1-d_2=1$.
As before, it suffices to prove that if $\dim \mathbb C[V]^G_{|H|}=\dim \mathbb C[V']^H_{|G|}$, then we have $\varphi_G(d)=\varphi_H(d)$ holds for any $d|(|G|,|H|)$. Assume to the contrary that $\varphi_G(d)\neq\varphi_H(d)$ for some $d|(|G|,|H|)$. Let $|G|=n$ and $|H|=m$.

Let
$$a=\min\{d\ |\ \varphi_G(d)\neq \varphi_H(d)\text{ for }d|(n,m)\}.$$

{\bf Case 1:} Assume that $\varphi_G(a)< \varphi_H(a).$
 If $\varphi_G(d)\leq \varphi_H(d)$ holds for any $d|(n,m)$, then the desired result follows. Therefore, we assume that
 $\varphi_G(d)>\varphi_H(d)$ holds for some $d|(n,m)$ and
let $$b=\min\{d\ |\ \varphi_G(d)> \varphi_H(d)\text{ for }d|(n,m)\}.$$
Clearly, we have $a<b$.
Now, we show that
\begin{equation}\label{finitegroupineq1}
(\varphi_H(a)-\varphi_G(a))
\binom{\frac{n+m}{a}}{\frac{n}{a},\frac{m}{a}}
>\left(\sum_{e\ge b}\varphi_{G}(e)\right)
\binom{\frac{n+m}{b}}{\frac{n}{b},\frac{m}{b}}.
\end{equation}
In fact, \eqref{finitegroupineq1} follows from the following stronger result
\begin{equation}\label{finitegroupineq2}
\binom{\frac{n+m}{a}}{\frac{n}{a},\frac{m}{a}}
>n\binom{\frac{n+m}{b}}{\frac{n}{b},\frac{m}{b}}.
\end{equation}
In order to prove \eqref{finitegroupineq2}, by Lemma \ref{lemma:Binom_ineq}.(i), we only need to prove that
\begin{equation}\label{finitegroupineq3}
(1+\frac{m}{n})^{n(\frac{1}{a}-\frac{1}{b})}(1+\frac{a}{b}\frac{n}{m})^{m(\frac{1}{a}-\frac{1}{b})}>n.
\end{equation}
It is easy to see that, as $\frac{n}{a}-\frac{n}{b}\ge 1$ and $m\ge n^2$,
$$(1+\frac{m}{n})^{n(\frac{1}{a}-\frac{1}{b})}(1+\frac{a}{b}\frac{n}{m})^{m(\frac{1}{a}-\frac{1}{b})}>1+\frac{m}{n}>n,$$
and therefore \eqref{finitegroupineq2} follows. By \eqref{finitegroupineq1} and the formula \eqref{eq:dimFormula}, it is easy to see that
$\dim \mathbb C[V']^H_{|G|}>\dim \mathbb C[V]^G_{|H|}$, a contradiction.

{\bf Case 2:} Assume that $\varphi_G(a)> \varphi_H(a)$.
 If $\varphi_G(d)\geq \varphi_H(d)$ holds for any $d|(n,m)$, then the desired result follows. Therefore, we assume that
 $\varphi_G(d)<\varphi_H(d)$ holds for some $d|(n,m)$ and
let $$b=\min\{d\ |\ \varphi_G(d)< \varphi_H(d)\text{ for }d|(n,m)\}.$$
Clearly, we have $a<b$.
Now, we show that
\begin{equation}\label{finitegroupineq4}
(\varphi_G(a)-\varphi_H(a))
\binom{\frac{n+m}{a}}{\frac{n}{a},\frac{m}{a}}
>\left(\sum_{e\ge b}\varphi_{H}(e)\right)
\binom{\frac{n+m}{b}}{\frac{n}{b},\frac{m}{b}}.
\end{equation}
In fact, \eqref{finitegroupineq4} follows from the following stronger result
\begin{equation}\label{finitegroupineq5}
\binom{\frac{n+m}{a}}{\frac{n}{a},\frac{m}{a}}
>m\binom{\frac{n+m}{b}}{\frac{n}{b},\frac{m}{b}}.
\end{equation}
In order to prove \eqref{finitegroupineq5}, by Lemma \ref{lemma:Binom_ineq}.(i), we only need to prove that
\begin{equation}\label{finitegroupineq6}
(1+\frac{m}{n})^{n(\frac{1}{a}-\frac{1}{b})}(1+\frac{a}{b}\frac{n}{m})^{m(\frac{1}{a}-\frac{1}{b})}>m.
\end{equation}
Note that, as $n$ does not contain two divisors $d_1, d_2>1$ with $d_1-d_2=1$, we have $n(\frac{1}{a}-\frac{1}{b})\ge 2$.  Since $m\ge n^2$, we have that
$$\left(1+\frac{m}{n}\right)
^{n(\frac{1}{a}-\frac{1}{b})}\geq \left(1+\frac{m}{n}\right)^2=1+2\frac{m}{n}+\frac{m^2}{n^2}>m,$$
as desired. By \eqref{finitegroupineq4} and the formula \eqref{eq:dimFormula}, it is easy to see that
$\dim \mathbb C[V]^G_{|H|}>\dim \mathbb C[V']^H_{|G|}$, a contradiction. This completes the proof.
\qed

\section{Further discussions}

In factorization theory (which is closely related to zero-sum theory; see \cite{GH}), the famous and fascinating \emph{characterization problem} has a similar flavor with Theorem \ref{thm:MainTheorem}. We recall some basic definitions; see \cite{Gero} for detailed explanations. Let $H$ be a commutative and cancellative monoid. If an element $a \in H$ has a factorization $a = u_1 \cdot \ldots \cdot u_k$, where $k \in \mathbb N$ and $u_1, \ldots, u_k$ are irreducible elements of $H$, then $k$ is called a factorization length of $a$ and the set $\mathsf L (a) \subset \mathbb N$ of all possible factorization lengths is called the {\it set of lengths} of $a$. We denote by $\mathcal L (H) = \{ \mathsf L (a) \colon a \in H \}$ the {\it system of sets of lengths} of $H$.

Now, let $H$ be a Krull monoid with  finite class group $G$ such that every class contains a prime divisor. Then there is a transfer homomorphism $\boldsymbol \beta \colon H \to \mathcal B (G)$, where $\mathcal B (G)$ is the monoid of zero-sum sequences over $G$, which implies that $\mathcal L (H) = \mathcal L \big( \mathcal B (G) \big)$. In particular, this shows that the system of sets of lengths of $H$ depends only on the class group $G$. The associated inverse problem (known as the Characterization Problem) asks whether the system $\mathcal L \big( \mathcal B (G) \big)$ is characteristic for the group. Indeed, the Characterization Problem runs as follows:

\emph{Let $G$ be a finite abelian group with Davenport constant $\mathsf D (G) \ge 4$, and let $G'$ be an abelian group with $\mathcal L \big( \mathcal B (G) \big) = \mathcal L \big( \mathcal B (G') \big)$. Are $G$ and $G'$ isomorphic?
}

The answer is affirmative for groups $G$ of rank $\mathsf r (G) \le 2$ (see \cite{Ge-Sc19a}) and the standing conjecture is that the answer is affirmative for all finite abelian groups $G$; we refer to the survey \cite{Ge-Zh20a} for an overview of this problem.

This paper provides some results regarding zero-sum sequences over finite abelian groups and polynomial invariants of finite groups. Recently, the relationship between zero-sum theory (also factorization theory) and invariant theory is getting closer; see \cite{CDG,CDS,Domo,HanZh} for some recent studies. Based on Theorems \ref{thm:MainTheorem}, \ref{Mainthm2}, and \ref{Mainthm3}, it is natural to propose the following conjecture.

\begin{conj}
Let $G$ and $H$ be finite groups. Let $V$ (resp. $V'$) be the regular representation of $G$ (resp. $H$). Then we have $$\dim \mathbb C[V]^G_{|H|}=\dim \mathbb C[V']^H_{|G|}$$
if and only if $\varphi_G(d)=\varphi_H(d)$ for any $d|(|G|,|H|)$.
\end{conj}

\subsection*{Acknowledgments}

We sincerely thank Prof. Alfred Geroldinger for helpful comments on our manuscript. M.S. Li was supported by the National Science Foundation of China Grant No.11871295. H.B. Zhang was supported by the National Science Foundation of China Grant No.11901563.

\bibliographystyle{amsplain}

\end{document}